\documentclass[preprint,12pt]{elsarticle}
\usepackage{amsfonts}
\usepackage{amsmath}
\usepackage{amssymb}

\input{tcilatex}

\begin{document}

\begin{frontmatter}



\title{
Fuzzy games with a countable space of actions and applications to systems of  generalized quasi-variational  inequalities
}

\author{Monica Patriche}

\address{
University of Bucharest
Faculty of Mathematics and Computer Science
    
14 Academiei Street
   
 010014 Bucharest, 
Romania
    
monica.patriche@yahoo.com }

\begin{abstract}
In this paper, we introduce an abstract fuzzy economy (generalized fuzzy game) model with a
 countable space of actions and we study the existence of the fuzzy equilibrium. As applications, 
two types of results are obtained. The first ones concern the existence of the solutions for systems 
of generalized quasi-variational inequalities with random fuzzy mappings which we define. The last 
ones are new random fixed point theorems for correspondences with values in complete countable 
metric spaces.
\end{abstract}

\begin{keyword}
Abstract fuzzy economy, \
fuzzy equilibrium, \
incomplete information, \
random fixed point, \
random quasi-variational inequalities, \
random fuzzy mapping.\


\end{keyword}

\end{frontmatter}



\label{}





\bibliographystyle{elsarticle-num}
\bibliography{<your-bib-database>}

\begin{thebibliography}{99}
\bibitem{ali} C. D. Aliprantis, K. C. Border, \textit{Infinite Dimensional
Analysis: A Hitchhiker's Guide.} Springer-Verlag, Berlin, 1994.

\bibitem{arr} Arrow K. J., Debreu G. (1954), \textit{Existence of an
Equilibrium for a Competitive} \textit{Economy}. Econometrica, \textbf{22},
265-290.

\bibitem{bor} Borglin A. and Keiding H. (1976), \textit{Existence of
equilibrium actions and of equilibrium: A note on the 'new' existence
theorem.} J. Math. Econom. \textbf{3}, 313-316.

\bibitem{4} S.-S. Chang and K. K. Tan, \textit{Equilibria and maximal
elements of abstract fuzzy economies and qualitative fuzzy games, }Fuzzy
Sets and Systems \textbf{125} (2002), 389-399.

\bibitem{3} Debreu G. (1952), A social equilibrium existence theorem, Proc.
Natl. Acad. Sci. USA \textbf{38}, 886-903.

\bibitem{fan} Fan K. (1952), Fixed-point and minimax theorems in locally
convex topological linear spaces, Proc. Nat. Acad. Sci. U.S.A. \textbf{38},
121-126.

\bibitem{6} Gwinner J. and Raciti F. (2006), On a class of random
variational inequalities on random sets, Numerical Functional Analysis and
Optimization, 27(5--6):619--636.

\bibitem{7} Gwinner J. and Raciti F. (2012), Some equilibrium problems under
uncertainty and random variational inequalities, Ann. Oper. Res., 200:
299-319.

\bibitem{har} Harsanyi, J.C.: Games with incomplete information played by
\textquotedblleft Bayesian\textquotedblright\ players. Management Science
14, 159--189 (1967).

\bibitem{8} Huang N.-J. (2001), \textit{Existence of equilibrium for
generalized abstract fuzzy economies,} Fuzzy Sets and Systems \textbf{117},
151-156.

\bibitem{11} Kim W. K. and Lee K. H. (1998), \textit{Fuzzy fixed point and
existence of equilibria of fuzzy games,} J. Fuzzy Math. \textbf{6}, 193-202.

\bibitem{12} Kim W. K. and Lee K. H.(1999), \textit{Existence of equilibria
in generalized fuzzy games, }J. Chungcheong Math. Soc. \textbf{12}, 53-61.

\bibitem{13} Kim W. K. and Lee K. H.(2001), \textit{Generalized fuzzy games
and fuzzy equilibria,} Fuzzy Sets and Systems \textbf{122}, 293-301.

\bibitem{nash} Nash J. F. (1950), \textit{Equilibrium points in n-person
games.} Proc. Nat. Acad. Sci. U. S. A., \textbf{36}, \textit{1}, 48-49.

\bibitem{noor} Noor M.A. and Elsanousi S.A.(1993), Iterative algorithms for
random variational inequalities, Pan Amer. Math. Journal 3, 39-50.

\bibitem{pat1} Patriche M. (2009), \textit{Equilibria for free abstract
fuzzy economies,} An. St. Ovidius University of Constanta, Volume 17(2),
143-154.

\bibitem{pat} Patriche M., Bayesian abstract economy with a measure space of
agents, Abstract and Applied Analysis, Volume 2009 (2009), Article ID 523619.

\bibitem{pat2} Patriche M.(2010), \textit{Abstract fuzzy economies and fuzzy
equilibrium pairs,} Mathematical Reports, 12(62), 4, 373 -- 382.

\bibitem{patt} Patriche M., Bayesian abstract economies, Optimization, Vol.
60, Nr 4, 2011, 451-461.

\bibitem{pat3} Patriche M. (2011a), \textit{Existence of fuzzy equilibria
for fuzzy abstract economies with Q'-majorized correspondences,} Recent
Researches in Neural Networks, Fuzzy Systems, Evolutionary Computing and
Automation, International Conference on Fuzzy Systems, 11 -- 13 April,
Brasov.

\bibitem{pat4} Patriche M. (2011b), \textit{Equilibrium in games and
competitive economies,} The Publishing House of the Romanian Academy.

\bibitem{pat5} Patriche M., Existence of equilibrium for an abstract economy
with private information and a countable space of actions, Mathematical
reports, to appear.

\bibitem{rad} R. Radner, 1968. Competitive Equilibrium Under Uncertainty,
Econometrica, 36 1, 31-58.

\bibitem{22} Shafer W. and Sonnenschein H. (1975), \textit{Equilibrium in
abstract economies without ordered} \textit{preferences, }Journal of
Mathematical Economics \textbf{2}, 345-348.

\bibitem{wang} Wang L., ChoY Je, Huang N. (2011), \textit{The robustness of
generalized abstract fuzzy economies in generalized convex spaces,} Fuzzy
Sets and Systems, 176 (1), 56-63.

\bibitem{19} Tan K.-K. and Yuan G. X.-Z. (1995), Random equilibria of random
generalized games with applications to non-compact random quasi-variational
inequalities, Topological Methods in Nonlinear Analysis, Journal of the
Juliusz Schauder Center, Volume 5, 59--82.

\bibitem{wang} L Wang, Y Je Cho, N Huang, \textit{The robustness of
generalized abstract fuzzy economies in generalized convex spaces,} Fuzzy
Sets and Systems, 176 (1) (2011), 56-63.

\bibitem{22} Yannelis N. C. and Prabhakar N. D. (1983), Existence of maximal
elements and equilibrium in linear topological spaces, J. Math. Econom. 
\textbf{12}, 233-245.

\bibitem{yu} H. Yu and Z. Zhang, \textit{Pure strategy equilibria in games
with countable actions,} Journal of Mathematical Economics 43 (2007),
192-200.

\bibitem{yuan} Yuan X. Z. (1988), The Study of Minimax inequalities and
Applications to Economies and Variational inequalities,\textit{\ }Memoirs of
the American Society, vol \textbf{132}, Nr 625.

\bibitem{z} L. A. Zadeh, \textit{Fuzzy sets,} Inform and Control \textbf{8}
(1965), 338-353.
\end{thebibliography}







\section{\textbf{INTRODUCTION}}

In equilibrium theory, there are two directions for the study of the models
with differential information. The first one is due to Harsanyi [9], whose
first approach refers to Nash type models, and the second one is due to
Radner [23], who introduced differential information into the Arrow-Debreu
model [2]. Since then, an entire literature has developed on this topic.

A basic problem concerning the noncooperative games with differential
information was the equilibrium existence. The classical deterministic
models of abstract economies due to Borglin and Keiding [3], Shafer and
Sonnenschein [24], or Yannelis and Prabhakar [28] have had several
generalizations. We can refer the reader to Yannelis's results, who
introduced new information concepts based on measurability requirements or
Patriche's models [17], [19], [21], [22].

A new approach of an abstract economy with private information and a
countable set of actions [22] follows the ideas of Yu and Zhang in [29], who
worked with the distributions of the atomless correspondences.

This paper has two aims. Firstly, we extend the study of the model defined
in [22] in the fuzzy framework, by taking into account those situations in
which agents may have partial control over the actions they choose. Since we
are interested in a fuzzy setting, special attention is paid to the fuzzy
equilibrium existence of the fuzzy games emphasizing that this approach can
be used to treat more complicated economic systems. The uncertainties
derived by the individual character of agents in election situations can be
interpreted using fuzzy random mappings, which can be used for the purpose
of modelling and analyzing the economic systems. We continue the way opened
by Kim and Lee [11], who proved that the theory of fuzzy sets, initiated by
Zadeh [31], has became a good framework for obtaining results concerning
fuzzy equlibrium existence for abstract fuzzy economies. For other results
in this domain, the reader is referred to [4], [10], [12], [13], [25], [27].

Finally, we search for applications of our results concerning the existence
of the fuzzy equilibrium for the fuzzy games with countable action sets, in
connection with the systems of generalized quasi-equilibrium inequalities
with random fuzzy mappings. The variational inequalities were introduced in
1960s by Fichera and Stampacchia, who studied equilibrium problems arising
from mechanics. Since then, this domain has been extensively developed and
has been found very useful in many diverse fields of pure and applied
sciences, such as mechanics, physics, optimization and control theory,
operations research and several branches of engineering sciences. Recently,
the random variational inequalities have been introduced and studied in [7],
[8], [15], [26], [30]. We prove the existence of the solutions for systems
of random generalized quasi-variational inequalities with random fuzzy
mappings which we define. As consequences, we obtain new random fixed point
theorems for correspondences with values in complete countable metric spaces.

The paper is organized as follows. In the next section, some notational and
terminological conventions are given. In Section 3, the model of an abstract
fuzzy economy with private information and a countable space of actions is
introduced and the main result is stated. Section 4 contains existence
theorems for solutions of random quasi-variational inequalities with random
fuzzy mappings and random fixed point theorems. In the last section,
classical results are obtained as consequences.

\section{\textbf{NOTATION AND DEFINITION}}

Throughout this paper, we shall use the following notation:

$%
\mathbb{R}
_{++}$ denotes the set of strictly positive reals. co$D$ denotes the convex
hull of the set $D$. $\overline{co}D$ denotes the closed convex hull of the
set $D$. $2^{D}$ denotes the set of all non-empty subsets of the set $D$. If 
$D\subset Y$, where $Y$ is a topological space, cl$D$ denotes the closure of 
$D$.\smallskip 

For the reader's convenience, we review a few basic definitions and results
concerning continuity and measurability of correspondences.

Let $X$, $Y$ be topological spaces and $F:X\rightarrow 2^{Y}$ be a
correspondence. $F$ is said to be \textit{upper semicontinuous} if for each $%
x\in X$ and each open set $V$ in $Y$ with $F(x)\subset V$, there exists an
open neighbourhood $U$ of $x$ in $X$ such that $F(y)\subset V$ for each $%
y\in U$. $F$ is said to be \textit{lower semicontinuous} \textit{(l.s.c)} if
for each $x\in X$ and each open set $V$ in $Y$ with $F(x)\cap V\neq
\emptyset $, there exists an open neighbourhood $U$ of $x$ in $X$ such that $%
F(y)\cap V\neq \emptyset $ for each $y\in U$.

Each correspondence\textbf{\ }$F:$ $X\rightarrow 2^{Y}$ has two natural
inverses. The upper inverse $F^{u}$ (also called the \textit{strong inverse}%
) of a subset $A$ of $Y$ is defined by $F^{u}(A)=\left\{ x\in A:F(x)\subset
A\right\} .$ The lower inverse $F^{l}$ (also called the \textit{weak inverse}%
) of a subset $A$ of $Y$ is defined by $F^{l}(A)=\left\{ x\in A:F(x)\cap
A\not=\emptyset \right\} .\medskip $

An important result concerning upper (lower) semicontinuous correspondences
is the following one.

\begin{lemma}
\cite{yuan} \textit{Let }$X$\textit{\ and }$Y$\textit{\ be two topological
spaces and let }$A$\textit{\ be a closed (resp. open) subset of }$X.$\textit{%
\ Suppose }$F_{1}:X\rightarrow 2^{Y}$\textit{, }$F_{2}:X\rightarrow 2^{Y}$%
\textit{\ are lower semicontinuous (resp. upper semicontinuous) such that }$%
F_{2}(x)\subset F_{1}(x)$\textit{\ for all }$x\in A.$\textit{\ Then the
correspondence }$F:X\rightarrow 2^{Y}$\textit{\ defined by}
\end{lemma}

\begin{center}
$\mathit{F(x)=}\left\{ 
\begin{array}{c}
F_{1}(x)\text{, \ \ \ \ \ \ \ if }x\notin A\text{, } \\ 
F_{2}(x)\text{, \ \ \ \ \ \ \ \ \ if }x\in A%
\end{array}%
\right. $
\end{center}

\textit{is also lower semicontinuous (resp. upper semicontinuous).\medskip }

Let $(T$, $\mathcal{T})$ be a measurable space, $Y$ a topological space and $%
F:T\rightarrow 2^{Y}$ a corespondence. $F$ is \textit{weakly} \textit{%
measurable }if $F^{l}(A)\in \mathcal{T}$ for each open subset $A$ of $Y.$ $F$
is \textit{measurable }if $F^{l}(A)\in \mathcal{T}$ for each closed subset $%
A $ of $Y.$ If $(T$, $\mathcal{T})$ is a measurable space, $Y$ a countable
set and $F:T\rightarrow 2^{Y}$ is a corespondence, then $F$ is measurable if
for each $y\in Y,$ $F^{-1}(y)=\left\{ t\in T:y\in F(t)\right\} $ is $%
\mathcal{T-} $measurable.\medskip

\begin{lemma}
\cite{ali}\textbf{. }For\textbf{\ }a correspondence $F:T\rightarrow 2^{Y}$
from a measurable space into a metrizable space we have the following:
\end{lemma}

\begin{enumerate}
\item If $F$ is measurable, then it is also weakly measurable$;$

\item If $F$ is compact valued and weakly measurable, it is measurable.$%
\medskip $
\end{enumerate}

The following properties are essential tools used to prove the existence of
equilibria for abstract economies in Section 3. We follow Yu and Zhang \cite%
{yu}. Let\textbf{\ }$Y$ be a countable complete metric space, $(T$, $%
\mathcal{T}$, $\lambda )$ an atomless probability space and $F:T\rightarrow
2^{Y}$ a measurable corespondence. The function $f:T\rightarrow Y$ is said
to be \textit{a selection of }$F$ if $f(t)\in F(t)$ for $\lambda -$almost $%
t\in T.$ Let us denote $\mathcal{D}_{F}=\left\{ \lambda f^{-1}:f\text{ is a
measurable selection of }F\right\} .\medskip $

We will present some regular properties of $\mathcal{D}_{F},$ also obtained
by Yu and Zhang \cite{yu}.

The next lemma states the convexity of $\mathcal{D}_{F}$ for any
correspondence $F$.

\begin{lemma}
\cite{yu}\textbf{. }Let\textbf{\ }$Y$ be a countable complete metric space, $%
(T$, $\mathcal{T}$, $\lambda )$ an atomless probability space and $%
F:T\rightarrow 2^{Y}$ a measurable corespondence. Then, $\mathcal{D}_{F}$ is
non-empty and convex in the space $\mathcal{M}(Y)$ - the space of
probability measure on $Y$, equipped with the topology of weak
convergence.\medskip 
\end{lemma}

The compactness of $\mathcal{D}_{F}$ is stated in the next result.

\begin{lemma}
\cite{yu}. Let\textbf{\ }$Y$ be a countable complete metric space, and $(T$, 
$\mathcal{T}$, $\lambda )$ be an atomless probability space and $%
F:T\rightarrow 2^{Y}$ be a measurable corespondence. If $F$ is compact
valued, then, $\mathcal{D}_{F}$ is compact in $\mathcal{M}(Y).\medskip $
\end{lemma}

Lemma 5 concerns the upper semicontinuity property of $\mathcal{D}_{F_{x}}.$

\begin{lemma}
\cite{yu}\textbf{.} Let $X$ be a metric space, $(T$, $\mathcal{T}$, $\lambda
)$ an atomless probability space, $Y$ a countable complete metric space and $%
F:T\times X\rightarrow 2^{Y}$ a correspondence. Let us assume that, for any
fixed $x$ in $X$, $F(\cdot ,x)$ (also denoted by $F_{x})$ is a
compact-valued measurable correspondence$,$ and for each fixed $t\in T,$ $%
F(t,\cdot )$ is upper semicontinuous on $X$. Also, let us assume that there
exists a compact valued corespondence $H:T\times X\rightarrow 2^{Y}$ such
that $F(t,x)\subset H(t)$ for all $t$ and $x$. Then $\mathcal{D}_{F_{x}}$ is
upper semicontinuous on $X.$
\end{lemma}

We will also need Kuratowski-Ryll-Nardzewski Selection Theorem in order to
prove our main results in the next section.

\begin{theorem}
(Kuratowski-Ryll-Nardzewski Selection \ Theorem)\textbf{\ \cite{ali}. }A
weakly measurable correspondence with non-empty closed values from a
measurable space into a Polish space admits a measurable selector.\medskip 
\end{theorem}

The fuzzy mappings are presented below.

Let $\tciFourier (Y)$ be a collection of all fuzzy sets over $Y.$ A mapping $%
F:\Omega \rightarrow \tciFourier (Y)$ is called a \textit{fuzzy mapping.} If 
$F$ is a fuzzy mapping from $\Omega ,$ $F(\omega )$ is a fuzzy set on $Y$
and $F(\omega )(y)$ is the membership function of $y$ in $F(\omega ).$

Let $A\in \tciFourier (Y),$ $a\in \lbrack 0,1],$ then the set $%
(A)_{a}=\{y\in Y:A(y)\geq a\}$ is called an $a-$cut set of the fuzzy set $A.$

A fuzzy mapping $F:\Omega \rightarrow \tciFourier (Y)$ is said to be \textit{%
measurable} if, for any given $a\in \lbrack 0,1],$ $(F(\cdot ))_{a}:\Omega
\rightarrow 2^{Y}$ is a measurable correspondence. A fuzzy mapping $F:\Omega
\rightarrow \tciFourier (Y)$ is said \textit{to have a measurable graph} if,
for any given $a\in \lbrack 0,1],$ the correspondence $(F(\cdot
))_{a}:\Omega \rightarrow 2^{Y}$ has a measurable graph. A fuzzy mapping $%
F:\Omega \times X\rightarrow \tciFourier (Y)$ is called a \textit{random
fuzzy mapping} if, for any given $x\in X,$ $F(\cdot ,x):\Omega \rightarrow
\tciFourier (Y)$ is a measurable fuzzy mapping.

\section{FUZZY EQUILIBRIUM EXISTENCE FOR ABSTRACT FUZZY ECONOMIES WITH
PRIVATE INFORMATION AND A COUNTABLE SET OF ACTIONS}

\subsection{THE MODEL OF AN ABSTRACT FUZZY ECONOMY WITH PRIVATE INFORMATION}

In this section we define a model of abstract fuzzy economy with private
information and a countable set of actions. This model is a generalization
of the one defined in \cite{pat5}

Let $I$ be a non-empty finite set (the set of agents). For each $i\in I$,
the space of actions, $S_{i}$ is a countable complete metric space and $%
(\Omega _{i},\mathcal{Z}_{i})$ is a measurable space. Let $(\Omega ,\mathcal{%
F})$ be the product measurable space $\underset{i\in I}{(\prod }\Omega _{i},%
\underset{i\in I}{\tbigotimes }\mathcal{Z}_{i})$, and let $\mu $ be a
probability measure on $(\Omega ,\mathcal{F}).$ For a point $\omega =(\omega
_{1},...,\omega _{n})\in \Omega ,$ define the coordinate projections $\tau
_{i}(\omega )=\omega _{i}.$ The random mapping $\tau _{i}(\omega )$ is
interpreted as player i's private information related to his action.

We also denote, for each $i\in I,$ Meas($\Omega _{i},S_{i})$ the set of
measurable fuzzy mappings $f$ from $(\Omega _{i},\mathcal{Z}_{i})$ to $%
S_{i}. $ An element $g_{i}$ of Meas($\Omega _{i},S_{i})$ is called a \textit{%
pure strategy} for player $i.$ A \textit{pure strategy profile} $g$ is an
n-vector function $(g_{1},g_{2},...,g_{n})$ that specifies a pure strategy
for each player.

We suppose that there exists a fuzzy mapping $X_{i}:\Omega _{i}\rightarrow 
\mathcal{F(}S_{i})$ such that each agent $i$ can choose an action from $%
(X_{i}(\omega _{i}))_{z_{i}}\subset S_{i}$ for each $\omega _{i}\in \Omega
_{i}.$

Let $D_{(X_{i}(\cdot ))_{z_{i}}}$\ be the set $\left\{ (\mu \tau
_{i}^{-1})g_{i}^{-1}:g_{i}\text{ is a measurable selection of }(X_{i}(\cdot
))_{z_{i}}\right\} $ and $\mathcal{D}_{(X(\cdot ))_{z}}:=\underset{i\in I}{%
\prod }\mathcal{D}_{(X_{i}(\cdot ))_{z_{i}}}.$

For each $i\in I,$ let us denote $h_{g_{i}}=(\mu \tau _{i}^{-1})g_{i}^{-1},$
where $g_{i}$ is a measurable selection of $(X_{i}(\cdot ))_{z_{i}}$ and $%
h_{g}=(h_{g_{1}},h_{g_{2}},...,h_{g_{n}}).$

\begin{definition}
An \textit{abstract fuzzy economy} \textit{(or a generalized fuzzy game)} 
\textit{with private information} \textit{and a countable space of actions}
is defined as
\end{definition}

$\Gamma =(I,((\Omega _{i},\mathcal{Z}%
_{i}),(S_{i},A_{i},P_{i},a_{i},b_{i},z_{i})))_{i\in I},\mu )),$ where:

(a) $X_{i}:\Omega _{i}\rightarrow \mathcal{F}(S_{i})$ is the \textit{action
(strategy) fuzzy mapping} of agent $i$,

(b) for each $\omega _{i}\in \Omega _{i},$ $A_{i}(\omega _{i},$%
\textperiodcentered $):\mathcal{D}_{(X(\cdot ))_{z}}\rightarrow \mathcal{F}%
(S_{i})$ is the \textit{random fuzzy constraint mapping of agent }$i;$

(c) for each $\omega _{i}\in \Omega _{i},$ $P_{i}(\omega _{i},$%
\textperiodcentered $):\mathcal{D}_{(X(\cdot ))_{z}}\rightarrow \mathcal{F}%
(S_{i})$ is the \textit{random fuzzy preference mapping of agent }$i;$

(d) $a_{i}:\mathcal{D}_{(X(\cdot ))_{z}}\rightarrow (0,1]$ is a \textit{%
random fuzzy constraint function} and $p_{i}:\mathcal{D}_{(X(\cdot
))_{z}}\rightarrow (0,1]$ is a \textit{random fuzzy preference function;}

(e) $z_{i}\in (0,1]$ is such that for all $(\omega _{i},h_{g})\in \Omega
_{i}\times \mathcal{D}_{(X(\cdot ))_{z}},$ $(A_{i}(\omega
_{i},h_{g}))_{a_{i}(h_{g})}\subset (X_{i}(\omega _{i}))_{z_{i}}$ and $%
(P_{i}(\omega ,h_{g}))_{p_{i}(h_{g})}\subset (X_{i}(\omega
))_{z_{i}}.\medskip $

\begin{definition}
A \textit{fuzzy} \textit{equilibrium} for $\Gamma $ is defined as a strategy
profile $g^{\ast }=(g_{1}^{\ast },g_{2}^{\ast },...,g_{n}^{\ast })\in
\tprod\limits_{i\in I}$Meas$(\Omega _{i},S_{i})$ such that for each $i\in I:$
\end{definition}

1) $g_{i}^{\ast }(\omega _{i})\in (A_{i}(\omega _{i},h_{g^{\ast
}})_{a_{i}(h_{g^{\ast }})}$ for each $\omega _{i}\in \Omega _{i};$

2) $(A_{i}(\omega _{i},h_{g^{\ast }}))_{a_{i}(h_{g^{\ast }})}\cap
(P_{i}(\omega _{i},h_{g^{\ast }}))_{p_{i}(h_{g^{\ast }})}=\phi $ for each $%
\omega _{i}\in \Omega _{i}.$

\subsection{EXISTENCE OF THE FUZZY EQUILIBRIUM FOR ABSTRACT FUZZY ECONOMIES
WITH A COUNTABLE SET OF ACTIONS}

In this subsection we prove the existence of fuzzy equilibrium of abstract
fuzzy economies.

Theorem 7 is our main result.

Let us denote $\lambda =(\lambda _{1},\lambda _{2},...,\lambda _{n})\in 
\mathcal{D}_{(X(\cdot ))_{z}}.$\medskip

\begin{theorem}
\textit{Let\ }$\Gamma =(I,((\Omega _{i},\mathcal{Z}%
_{i}),(S_{i},A_{i},P_{i},a_{i},b_{i},z_{i}))_{i\in I},\mu )$\textit{\ be an
abstract economy with private information and a countable space of actions,
where }$I$\textit{\ is a finite index set such that for each }$i\in I:$
\end{theorem}

\textit{a) }$S_{i}$\textit{\ is a countable complete metric space\ and }$%
(\Omega _{i},\mathcal{Z}_{i})$\textit{\ is a measurable space; }$(\Omega ,%
\mathcal{Z})$\textit{\ is the product measurable space }$\underset{i\in I}{%
(\prod }(\Omega _{i},\mathcal{Z}_{i}))$\textit{\ and }$\mu $\textit{\ an
atomless probability measure on }$(\Omega ,Z);$

\textit{b) the fuzzy mapping }$X_{i}:\Omega _{i}\rightarrow \mathcal{F}%
(S_{i})$\textit{\ is measurable and }$(X_{i}(\cdot ))_{z_{i}}:\Omega
_{i}\rightarrow 2^{S_{i}}$ \textit{has} \textit{compact values; }

\textit{c) the fuzzy mapping }$A_{i}$ \textit{is such that} \textit{for each}
$\lambda \in \mathcal{D}_{(X(\cdot ))_{z}},$ $\mathit{(}A_{i}(\cdot ,\lambda
)_{a_{i}(\lambda )}:\Omega _{i}\rightarrow 2^{S_{i}}$\textit{\ is
measurable\ and, for all }$\omega _{i}\in \Omega _{i},$\textit{\ }$%
(A_{i}(\omega _{i},\cdot ))_{a_{i}(\cdot )}:\mathcal{D}_{(X(\cdot
))_{z}}\rightarrow 2^{S_{i}}$\textit{\ is upper semicontinuous with
non-empty compact values;}

\textit{d) the fuzzy mapping }$P_{i}$ \textit{is such that} \textit{for each}
$\lambda \in \mathcal{D}_{(X(\cdot ))_{z}},$ $\mathit{(}P_{i}(\cdot ,\lambda
)_{p_{i}(\lambda )}:\Omega _{i}\rightarrow 2^{S_{i}}$\textit{\ is
measurable\ and, for all }$\omega _{i}\in \Omega _{i},$\textit{\ }$%
(P_{i}(\omega _{i},\cdot ))_{p_{i}(\cdot )}:\mathcal{D}_{(X(\cdot
))_{z}}\rightarrow 2^{S_{i}}$\textit{\ is upper semicontinuous with
non-empty compact values;}

\textit{e) for each }$\omega _{i}\in \Omega _{i}$\textit{\ and each }$g\in
\tprod\limits_{i\in I}$\textit{Meas}$(\Omega _{i},S_{i}),$ \textit{\ }$%
g_{i}(\omega _{i})\not\in (P_{i}(\omega _{i},h_{g})_{p_{i}((h_{g})};$

f) the set $U_{i}^{\omega _{i}}:=\left\{ \lambda \in \mathcal{D}_{(X(\cdot
))_{z}}:(A_{i}(\omega _{i},\lambda ))_{a_{i}(\lambda )}\cap P_{i}((\omega
_{i},\lambda ))_{p_{i}(\lambda )}=\emptyset \right\} $ is open in $\mathcal{D%
}_{(X(\cdot ))_{z}}$ in for each $\omega _{i}\in \Omega _{i}$.

\textit{Then, there exists }$g^{\ast }\in \tprod\limits_{i\in I}$\textit{Meas%
}$(\Omega _{i},S_{i})$\textit{\ an equilibrium for }$\Gamma .$\textit{%
\bigskip }

\textit{Proof.} According to Lemma 3, $D_{(X_{i}(\cdot ))_{z_{i}}}$ is
non-empty and convex. According to Lemma 4, $D_{(X_{i}(\cdot ))_{z_{i}}}$ is
compact$.$ For each $i\in I$ let us denote

$U_{i}:=\left\{ (\omega _{i},\lambda )\in \Omega _{i}\times \mathcal{D}%
_{(X(\cdot ))_{z}}:(A_{i}(\omega _{i},\lambda ))_{a_{i}(\lambda )}\cap
P_{i}((\omega _{i},\lambda ))_{p_{i}(\lambda )}=\emptyset \right\} $ and

$U_{i}^{\omega _{i}}:=\left\{ (\omega _{i},\lambda )\in \Omega _{i}\times 
\mathcal{D}_{(X(\cdot ))_{z}}:(A_{i}(\omega _{i},\lambda ))_{a_{i}(\lambda
)}\cap P_{i}((\omega _{i},\lambda ))_{p_{i}(\lambda )}=\emptyset \right\} .$

According to the assumption f), the set $U_{i}^{\omega _{i}}$ is open in $%
\mathcal{D}_{(X(\cdot ))_{z}}.$

Let us define $F_{i}:\Omega _{i}\times \mathcal{D}_{(X(\cdot
))_{z}}\rightarrow 2^{S_{i}}$ by

$F_{i}(\omega _{i},\lambda )=\left\{ 
\begin{array}{c}
(A_{i}(\omega _{i},\lambda ))_{a_{i}(\lambda )}\cap (P_{i}(\omega
_{i},\lambda ))_{p_{i}(\lambda )}\text{ if }(\omega _{i},\lambda )\not\in
U_{i}, \\ 
(A_{i}(\omega _{i},\lambda ))_{a_{i}(\lambda )}\text{ if }(\omega
_{i},\lambda )\in U_{i}.%
\end{array}%
\right. $

Then, the correspondence $F_{i}$ has non-empty compact values and is
measurable with respect to $\Omega _{i}$ and upper semicontinuous with
respect to\textit{\ }$\lambda \in \mathcal{D}_{(X(\cdot ))_{z}}.$

We denote $\mathcal{D}_{F_{i}}(\lambda )=$

=$\{h_{g_{i}}=(\mu \tau _{i}^{-1})g_{i}^{-1}:g_{i}$ is a measurable
selection of $F_{i}(\cdot ,\lambda )\}.$ Then:

i) $\mathcal{D}_{F_{i}}(\lambda )$ is non-empty because there exists a
measurable selection from the correspondence $F_{i}$ according to
Kuratowski-Ryll-Nardewski Selection Theorem.

ii) $\ \mathcal{D}_{F_{i}}(\lambda )$ is convex and compact according to
Lemma 3 and Lemma 4.

We define $\Phi :$ $\mathcal{D}_{(X(\cdot ))_{z}}\rightarrow 2^{\mathcal{D}%
_{(X(\cdot ))_{z}}},$ $\Phi (\lambda )=\tprod\limits_{i\in I}\mathcal{D}%
_{F_{i}}(\lambda ).$

The set $\mathcal{D}_{(X(\cdot ))_{z}}$ is non-empty compact and convex.
According to Lemma 5, the correspondence $\mathcal{D}_{F_{i}}$ is upper
semicontinuous. Then, the correspondence $\Phi $ is upper semicontinuous and
has non-empty compact and convex values. According to Ky Fan fixed point
Theorem [6], there exists a fixed point $\lambda ^{\ast }\in \Phi (\lambda
^{\ast }).$ In particular, for each player $i,$ $\lambda _{i}^{\ast }\in 
\mathcal{D}_{F_{i}}(\lambda ^{\ast }).$ Therefore, for each player $i,$
there exists $g_{i}^{\ast }\in $\textit{Meas}$(\Omega _{i},S_{i})$ such that 
$g_{i}^{\ast }$ is a selection of $F_{i}(\cdot ,\lambda ^{\ast })$ and $%
h_{g_{i}^{\ast }}=(\mu \tau _{i}^{-1})(g_{i}^{\ast })^{-1}=\lambda
_{i}^{\ast }.$ Let us denote $h_{g^{\ast }}=(h_{g_{1}^{\ast
}},...,h_{g_{n}^{\ast }}).$

We prove that $g^{\ast }$ is an equilibrium for $\Gamma .$ For each $i\in I,$
because $g_{i}^{\ast }$ is a selection of $F_{i}(\cdot ,h_{g_{1}^{\ast
}},...,h_{g_{n}^{\ast }})$, it follows that $g_{i}^{\ast }(\omega _{i})\in
(A_{i}(\omega _{i},h_{g^{\ast }}))_{a_{i}(h_{g^{\ast }})}\cap (P_{i}(\omega
_{i},h_{g^{\ast }}))_{p_{i}(h_{g^{\ast }})}$ if $(\omega _{i},h_{g^{\ast
}})\not\in U_{i}$ or $g_{i}^{\ast }(\omega _{i})\in (A_{i}(\omega
_{i},h_{g^{\ast }}))_{a_{i}(h_{g^{\ast }})}$ if $(\omega _{i},h_{g^{\ast
}})\in U_{i}.$

According to the assumption d), it follows that $g_{i}^{\ast }(\omega
_{i})\not\in (P_{i}(\omega _{i},h_{g^{\ast }}))_{p_{i}(h_{g^{\ast }})}$ for
each $\omega _{i}\in \Omega _{i}.$ Then $g_{i}^{\ast }(\omega _{i})\in
(A_{i}(\omega _{i},h_{g^{\ast }}))_{a_{i}(h_{g^{\ast }})}$ and $(\omega
_{i},h_{g^{\ast }})\in U_{i}.$ This is equivalent with the fact that $%
g_{i}^{\ast }(\omega _{i})\in (A_{i}(\omega _{i},h_{g^{\ast
}})_{a_{i}(h_{g^{\ast }})}$ and $(A_{i}(\omega _{i},h_{g^{\ast
}})_{a_{i}(h_{g^{\ast }})}\cap P_{i}(\omega _{i},h_{g^{\ast
}})_{p_{i}(h_{g^{\ast }})}=\emptyset $ for each $\omega _{i}\in \Omega _{i}.$
Consequently, $g^{\ast }=(g_{1}^{\ast },g_{2}^{\ast },...,g_{n}^{\ast })$ is
an equilibrium for $\Gamma .\medskip $

\section{ RANDOM\ QUASI-VARIATIONAL INEQUALITIES WITH RANDOM FUZZY MAPPINGS}

\subsection{New types of systems of generalized quasi-variational
inequalities}

Noor and Elsanousi \cite{noor} introduced the notion of a random variational
inequality. We also propose the next systems of generalized
quasi-variational inequalities.

Let $I$ be a non-empty and finite set. For each $i\in I$, let $S_{i}$ be a
countable complete metric space and $(\Omega _{i},\mathcal{Z}_{i})$ a
measurable space. Let $(\Omega ,\mathcal{F})$ be the product measurable
space $\underset{i\in I}{(\prod }\Omega _{i},\underset{i\in I}{\tbigotimes }%
\mathcal{Z}_{i})$ and $\mu $ a probability measure on $(\Omega ,\mathcal{F}).
$ For a point $\omega =(\omega _{1},...,\omega _{n})\in \Omega ,$ define the
coordinate projections $\tau _{i}(\omega )=\omega _{i}.$We also denote, for
each $i\in I,$ Meas($\Omega _{i},S_{i})$ the set of measurable fuzzy
mappings $f$ from $(\Omega _{i},\mathcal{Z}_{i})$ to $S_{i}.$ We suppose
that there exists a correspondence $X_{i}:\Omega _{i}\rightarrow S_{i}$ and
let $\mathcal{D}_{X_{i}}$\ be the set \{$(\mu \tau _{i}^{-1})g_{i}^{-1}:g_{i}
$ is a measurable selection of $X_{i}\}$ and $\mathcal{D}_{X}:=\underset{%
i\in I}{\prod }\mathcal{D}_{X_{i}}.$

Let $A_{i}:\Omega _{i}\times \mathcal{D}_{X}\rightarrow 2^{S_{i}}$ be a
correspondence and $\psi _{i}:\Omega _{i}\times \mathcal{D}_{X}\times
S_{i}\rightarrow \mathbb{R}\cup \{-\infty ,+\infty \}$.

We associate the next generalized quasi-variational problem with $A_{i}$ and 
$\psi _{i}$:\medskip 

(1) Find $\lambda ^{\ast }\in \mathcal{D}_{X}$ such that:

i) $\lambda _{i}^{\ast }(\omega _{i})\in A_{i}(\omega _{i},\lambda ^{\ast
}); $

ii) $\sup_{y_{i}\in A_{i}(\omega _{i},\lambda ^{\ast })}\psi _{i}(\omega
_{i},\lambda ^{\ast },y_{i})\leq 0$ for all $\omega _{i}\in \Omega _{i}.$%
\medskip

If $S_{i}$ is a countable completely metrizable topological vector space, we
introduce the following definition.

Let $S_{i}^{\prime }$ be the dual space of $S_{i}$ and $G_{i}:\Omega
_{i}\times S_{i}\rightarrow 2^{S_{i}^{\prime }},$ $A_{i}:\Omega _{i}\times 
\mathcal{D}_{X}\rightarrow 2^{S_{i}}$ be correspondences.

(2) Find $\lambda ^{\ast }\in \mathcal{D}_{X}$ such that:

i) $\lambda ^{\ast }(\omega _{i})\in A_{i}(\omega _{i},\lambda ^{\ast })$

ii) $\sup_{y_{i}\in A_{i}(\omega _{i},\lambda ^{\ast })}\sup_{v\in
G_{i}(\omega _{i},y_{i})}\func{Re}\langle v,\lambda _{i}^{\ast }(\omega
_{i})-y_{i}\rangle \leq 0$ for all $\omega _{i}\in \Omega _{i},$

where the real part of pairing between $S_{i}^{\prime }$ and $S_{i}$ is
denoted by $\func{Re}\langle v,x\rangle $ for each $v\in S_{i}^{\prime }$
and $x\in S_{i}.$\medskip

We will work on the next fuzzy model:\medskip

For each $i\in I$, let $S_{i}$ be a countable complete metric space and $%
(\Omega _{i},\mathcal{Z}_{i})$ is a measurable space. Let $A_{i}:\Omega
_{i}\times \mathcal{D}_{X}\rightarrow \mathcal{F}(S_{i})$ be a fuzzy mapping
and let $a_{i}:\mathcal{D}_{X}\rightarrow (0,1]$ be a fuzzy function. Let $%
\psi _{i}:\Omega _{i}\times \mathcal{D}_{X}\times S_{i}\rightarrow \mathbb{R}%
\cup \{-\infty ,+\infty \}$.

Now, we are introducing the next type of variational inequality:\medskip

(3) Find $\lambda ^{\ast }\in \mathcal{D}_{X}$ such that for every\textit{\ }%
$i\in I$\textit{,}

i) $\lambda _{i}^{\ast }(\omega _{i})\in (A_{i}(\omega _{i},\lambda ^{\ast
}))_{a_{i}(\lambda ^{\ast })};$

ii) $\sup_{y_{i}\in (A_{i}((\omega _{i},\lambda ^{\ast })))_{a_{i}(\lambda
^{\ast })}}\psi _{i}(\omega _{i},\lambda ^{\ast },y_{i})\leq 0$ for all $%
\omega _{i}\in \Omega _{i}.$

where $(A_{i_{x^{\ast }}})_{a_{i}(x^{\ast })}=\{z\in Y_{i}:A_{i_{x^{\ast
}}}(z)\geq a_{i}(x^{\ast })\}.$\medskip

If $A_{i}:X\rightarrow 2^{S_{i}}$ is a classical correspondence, then we get
the variational inequality defined in (1).\medskip 

Finally, we introduce the following system of generalized quasi-variational
inequalities, in case that, for each $i\in I,$ $S_{i}$ is a countable
completely metrizable topological vector space, $(\Omega _{i},\mathcal{Z}%
_{i})$ is measurable space, $A_{i}:\Omega _{i}\times \mathcal{D}%
_{X}\rightarrow \mathcal{F}(S_{i})$ and $G_{i}:\Omega _{i}\times
S_{i}\rightarrow \mathcal{F}(S_{i}^{\prime })$ are fuzzy mappings and $a_{i}:%
\mathcal{D}_{X}\rightarrow (0,1]$, $g_{i}:\mathcal{S}_{i}\rightarrow (0,1]$
are fuzzy functions.

(4) Find $\lambda ^{\ast }\in \mathcal{D}_{X}$ such that for every\textit{\ }%
$i\in I$\textit{,}

$\left\{ 
\begin{array}{c}
\lambda _{i}^{\ast }(\omega _{i})\in (A_{i}(\omega _{i},\lambda ^{\ast
}))_{a_{i}(\lambda ^{\ast })}; \\ 
\sup_{y_{i}\in (A_{i}(\omega _{i},\lambda ^{\ast }))_{a_{i}(\lambda ^{\ast
})}}\sup_{v\in (G_{i}(\omega _{i},y_{i}))_{g_{i}(y_{i})}}\func{Re}\langle
v,\lambda _{i}^{\ast }(\omega _{i})-y_{i}\rangle \leq 0\text{ }%
\end{array}%
\right. $ for all $\omega _{i}\in \Omega _{i}.$

\textit{\medskip }

\subsection{The existence of the solutions of the systems of generalized
quasi-variational inequalities with random fuzzy mappings}

In this section, we are establishing new results concerning the existence of
the systems of generalized random quasi-variational inequalities with random
fuzzy mappings and we also are stating random fixed point theorems. The
proofs rely on the theorem of fuzzy equilibrium existence for the abstract
fuzzy economy.

This is our first theorem.

\begin{theorem}
Let $I$ be a non-empty and finite set. For each $i\in I$, $S_{i}$ is a
countable complete metric space and $(\Omega _{i},\mathcal{Z}_{i})$ is a
measurable space. Let $(\Omega ,\mathcal{Z})$ be the product measurable
space $\underset{i\in I}{(\prod }\Omega _{i},\underset{i\in I}{\tbigotimes }%
\mathcal{Z}_{i})$, and let $\mu $ be a probability measure on $(\Omega ,%
\mathcal{Z}).$\textit{\ Suppose that the following conditions are satisfied
for each }$i\in I:$
\end{theorem}

\textit{a) the correspondence }$X_{i}:\Omega _{i}\rightarrow \mathcal{F}%
(S_{i})$\textit{\ is measurable such that }$(X_{i}(\cdot ))_{z_{i}}:\Omega
_{i}\rightarrow 2^{S_{i}}$ \textit{has compact values; }

\textit{b) the fuzzy mapping }$A_{i}$ \textit{is such that} \textit{for each}
$\lambda \in \mathcal{D}_{(X(\cdot ))_{z}},$ $\mathit{(}A_{i}(\cdot ,\lambda
)_{a_{i}(\lambda )}:\Omega _{i}\rightarrow 2^{S_{i}}$\textit{\ is
measurable\ and, for all }$\omega _{i}\in \Omega _{i},$\textit{\ }$%
(A_{i}(\omega _{i},\cdot ))_{a_{i}(\cdot )}:\mathcal{D}_{(X(\cdot
))_{z}}\rightarrow 2^{S_{i}}$\textit{\ is upper semicontinuous with
non-empty compact values;}

Let us \textit{assume that the mapping }$\psi _{i}:\Omega _{i}\times 
\mathcal{D}_{(X(\cdot ))_{z}}\times S_{i}\rightarrow \mathbb{R}\cup
\{-\infty ,+\infty \}$\textit{\ is such that:}

\ \ \ (\textit{c) }$\lambda \rightarrow \{y\in S_{i}:\psi _{i}(\omega
,\lambda ,y)>0\}:\mathcal{D}_{(X(\cdot ))_{z}}\rightarrow 2^{S_{i}}$\textit{%
\ is upper semicontinuous with compact values on }$\mathcal{D}_{(X(\cdot
))_{z}}$\textit{\ for each fixed }$\omega _{i}\in \Omega _{i};$

\ \ \ (\textit{d) }$\lambda _{i}(\omega _{i})\notin \{y\in S_{i}:\psi
_{i}(\omega _{i},\lambda ,y)>0\}$\textit{\ for each fixed }$(\omega
_{i},\lambda )\in \Omega _{i}\times \mathcal{D}_{(X(\cdot ))_{z}};$

\ \ \ (\textit{e) for each }$\omega _{i}\in \Omega _{i},$\textit{\ }$%
\{\lambda \in \mathcal{D}_{(X(\cdot ))_{z}}:\alpha _{i}(\omega _{i},\lambda
)>0\}$\textit{\ is weakly open in} $\mathcal{D}_{(X(\cdot ))_{z}},$\textit{\
where }$\alpha _{i}:\Omega _{i}\times \mathcal{D}_{(X(\cdot
))_{z}}\rightarrow \mathbb{R}$\textit{\ is defined by }

$\alpha _{i}(\omega _{i},\lambda )=\sup_{y\in (A_{i}(\omega _{i},\lambda
))_{a_{i}(\lambda ))}}\psi _{i}(\omega _{i},\lambda ,y)$\textit{\ for each }$%
(\omega _{i},\lambda )\in \Omega _{i}\times \mathcal{D}_{(X(\cdot ))_{z}};$

\ \ \ (\textit{f) }$\{\omega _{i}:\alpha _{i}(\omega _{i},\lambda )>0\}\in 
\mathcal{Z}_{i}$ for each $\lambda \in \mathcal{D}_{(X(\cdot ))_{z}}$\textit{%
.}

\textit{Then, there exists }$\lambda ^{\ast }\in \mathcal{D}_{(X(\cdot
))_{z}}$\textit{\ such that for every} $i\in I$,

i) $\lambda _{i}^{\ast }(\omega _{i})\in (A_{i}(\omega _{i},\lambda ^{\ast
}))_{a_{i}(\lambda ^{\ast })};$

\textit{ii) sup}$_{y\in (A_{i}(\omega ,\lambda ^{\ast })_{a_{i}(\lambda
^{\ast })}}\psi _{i}(\omega _{i},\lambda ^{\ast },y)\leq 0$\textit{.}$%
\medskip $

\textit{Proof.} For every $i\in I,$ let $P_{i}:\Omega _{i}\times \mathcal{D}%
_{(X(\cdot ))_{z}}\rightarrow \mathcal{F}(S_{i})$ and $p_{i}:\mathcal{D}%
_{(X(\cdot ))_{z}}\rightarrow (0,1]$ such that $(P_{i}(\omega ,\lambda
))_{p_{i}(\lambda )}=\{y\in S_{i}:\psi _{i}(\omega _{i},\lambda ,y)>0\}$ for
each $(\omega _{i},\lambda )\in \Omega _{i}\times \mathcal{D}_{(X(\cdot
))_{z}}.$

We shall show that the abstract economy $\Gamma =(I,((\Omega _{i},\mathcal{F}%
_{i}),(S_{i},A_{i},P_{i},a_{i},b_{i},$

\noindent $z_{i})))_{i\in I},\mu ))$ satisfies all the hypotheses of Theorem
7.

Suppose $\omega _{i}\in \Omega _{i}.$

According to c), we have that\textit{\ }$\lambda \rightarrow (P_{i}(\omega
_{i},\lambda ))_{p_{i}((\lambda )}:\mathcal{D}_{(X(\cdot ))_{z}}\rightarrow
2^{S_{i}}$\textit{\ }is upper semicontinuous\textit{\ }with non-empty values
and according to d), $\lambda _{i}(\omega _{i})\not\in (P_{i}(\omega
_{i},\lambda ))_{p_{i}((\lambda )}$ for each $\lambda \in \mathcal{D}%
_{(X(\cdot ))_{z}}.$

According to the definition of $\alpha _{i},$ we note that, for each $\omega
_{i}\in \Omega _{i},$ $\ \{\lambda \in \mathcal{D}_{(X(\cdot
))_{z}}:(A_{i}(\omega _{i},\lambda ))_{a_{i}(\lambda ))}\cap (P_{i}(\omega
_{i},\lambda ))_{p_{i}(\lambda )}\neq \emptyset \}=\{\lambda \in \mathcal{D}%
_{(X(\cdot ))_{z}}:\alpha _{i}(\omega ,\lambda )>0\}$ so that $\{\lambda \in 
\mathcal{D}_{(X(\cdot ))_{z}}:(A_{i}(\omega _{i},\lambda ))_{a_{i}(\lambda
)}\cap (P_{i}(\omega _{i},\lambda ))_{p_{i}(\lambda )}\neq \emptyset \}$ is
weakly open in $\mathcal{D}_{(X(\cdot ))_{z}}$ by e).

According to b) and f), it follows that for each $\lambda \in \mathcal{D}%
_{(X(\cdot ))_{z}},$ the correspondences\textit{\ }$(A_{i}(\cdot ,\lambda
))_{a_{i}(\lambda )}:\Omega _{i}\rightarrow 2^{S_{i}}$ and $(P_{i}(\omega
_{i},\lambda ))_{p_{i}(\lambda )}:\Omega _{i}\rightarrow 2^{S_{i}}$\textit{\ 
}are measurable$.$

Thus, the abstract fuzzy economy $\Gamma =(I,((\Omega _{i},\mathcal{Z}%
_{i}),(S_{i},A_{i},P_{i},a_{i},b_{i},z_{i})))_{i\in I},$

\noindent $\mu ))$ satisfies all the hypotheses of Theorem 7. Therefore,
there exists $\lambda ^{\ast }\in \mathcal{D}_{(X(\cdot ))_{z}}$ such that
for every $i\in I:$

$\lambda _{i}^{\ast }(\omega _{i})\in (A_{i}(\omega _{i},\lambda ^{\ast
}))_{a_{i}(\lambda ^{\ast })}$ and

$(A_{i}(\omega _{i},\lambda ^{\ast }))_{a_{i}(\lambda ^{\ast })}\cap
(P_{i}(\omega _{i},\lambda ^{\ast }))_{p_{i}(\lambda ^{\ast })}=\phi $;

that is, there exists $\lambda ^{\ast }\in \mathcal{D}_{(X(\cdot ))_{z}}$
such that for every $i\in I:$

i) $\lambda _{i}^{\ast }(\omega )\in (A_{i}(\omega ,\lambda ^{\ast
}))_{a_{i}\lambda ^{\ast })};$

ii) sup$_{y\in (A_{i}(\omega _{i},\lambda ^{\ast }))_{a_{i}(\lambda ^{\ast
})}}\psi _{i}(\omega _{i},\lambda ^{\ast },y)\leq 0$.$\medskip $

If \TEXTsymbol{\vert}I\TEXTsymbol{\vert}=1, we obtain the following
corollary.

\begin{corollary}
Let $S$ be a countable complete metric space, $(\Omega ,\mathcal{Z},\mathcal{%
\mu })$ be a measure space$.$\textit{\ Suppose that the following conditions
are satisfied}$:$
\end{corollary}

\textit{a) the fuzzy mapping }$X:\Omega \rightarrow \mathcal{F}(S)$\textit{\
is measurable such that }$(X(\cdot ))_{z}:\Omega _{i}\rightarrow 2^{S}$ 
\textit{has compact values;}

\textit{b) the fuzzy mapping }$A$ \textit{is such that} \textit{for each} $%
\lambda \in \mathcal{D}_{(X(\cdot ))_{z}},$ $\mathit{(}A(\cdot ,\lambda
))_{a(\lambda )}:\Omega \rightarrow 2^{S}$\textit{\ is measurable\ and, for
all }$\omega \in \Omega ,$\textit{\ }$(A(\omega ,\cdot ))_{a(\cdot )}:%
\mathcal{D}_{(X(\cdot ))_{z}}\rightarrow 2^{S}$\textit{\ is upper
semicontinuous with non-empty compact values;}

\textit{The mapping }$\psi :\Omega \times \mathcal{D}_{(X(\cdot
))_{z}}\times S\rightarrow \mathbb{R}\cup \{-\infty ,+\infty \}$\textit{\ is
such that:}

\ \ \ (\textit{c) }$\lambda \rightarrow \{y\in Y:\psi (\omega ,\lambda
,y)>0\}:\mathcal{D}_{(X(\cdot ))_{z}}\rightarrow 2^{S}$\textit{\ is upper
semicontinuous with compact values on }$\mathcal{D}_{(X(\cdot ))_{z}}$%
\textit{\ for each fixed }$\omega \in \Omega ;$

\ \ \ (\textit{d) }$\lambda (\omega )\notin \{y\in S:\psi (\omega ,\lambda
,y)>0\}$\textit{\ for each fixed }$(\omega ,\lambda )\in \Omega \times 
\mathcal{D}_{(X(\cdot ))_{z}};$

\ \ \ (\textit{e) for each }$\omega \in \Omega ,$\textit{\ }$\{\lambda \in 
\mathcal{D}_{(X(\cdot ))_{z}}:\alpha (\omega ,\lambda )>0\}$\textit{\ is
weakly open in} $\mathcal{D}_{(X(\cdot ))_{z}},$\textit{\ where }$\alpha
:Z\times \mathcal{D}_{(X(\cdot ))_{z}}\rightarrow \mathbb{R}$\textit{\ is
defined by }

$\alpha (\omega ,\lambda )=\sup_{y\in (A(\omega ,\lambda ))_{a(\lambda
))}}\psi (\omega ,\lambda ),y)$\textit{\ for each }$(\omega ,\lambda )\in
Z\times \mathcal{D}_{(X(\cdot ))_{z}};$

\ \ \ (\textit{f) }$\{\omega :\alpha (\omega ,\lambda )>0\}\in \mathcal{Z}$
for each $\lambda \in \mathcal{D}_{(X(\cdot ))_{z}}.$

\textit{Then, there exists }$\lambda ^{\ast }\in \mathcal{D}_{(X(\cdot
))_{z}}$\textit{\ such that},

i) $\lambda ^{\ast }(\omega )\in (A(\omega ,\lambda ^{\ast
}))_{a_{i}(\lambda ^{\ast })};$

\textit{ii) sup}$_{y\in (A(\omega ,\lambda ^{\ast })_{a(\lambda ^{\ast
})}}\psi (\omega ,\lambda ^{\ast },y)\leq 0$\textit{.}$\medskip $

As a consequence of Theorem 8, we prove the following Tan and Yuan's type
(1995) random quasi-variational inequality with random fuzzy mappings.

\begin{theorem}
Let $I$ be a non-empty and finite set. For each $i\in I$, $S_{i}$ is a
countable completely metrizable topological vector space and $(\Omega _{i},%
\mathcal{Z}_{i})$ is a measurable space. Let $(\Omega ,\mathcal{Z})$ be the
product measurable space $\underset{i\in I}{(\prod }\Omega _{i},\underset{%
i\in I}{\tbigotimes }\mathcal{Z}_{i})$, and $\mu $ a probability measure on $%
(\Omega ,\mathcal{Z}).$\textit{\ Suppose that the following conditions are
satisfied for each }$i\in I:$
\end{theorem}

\textit{a) the correspondence }$X_{i}:\Omega _{i}\rightarrow \mathcal{F}%
(S_{i})$\textit{\ is measurable such that }$(X_{i}(\cdot ))_{z_{i}}:\Omega
_{i}\rightarrow 2^{S_{i}}$ \textit{has compact values; }

\textit{b) the mapping }$A_{i}$ \textit{is such that} \textit{for each} $%
\lambda \in \mathcal{D}_{(X(\cdot ))_{z}},$ $\mathit{(}A_{i}(\cdot ,\lambda
)_{a_{i}(\lambda )}:\Omega _{i}\rightarrow 2^{S_{i}}$\textit{\ is
measurable\ and, for all }$\omega _{i}\in \Omega _{i},$\textit{\ }$%
(A_{i}(\omega _{i},\cdot ))_{a_{i}(\cdot )}:\mathcal{D}_{(X(\cdot
))_{z}}\rightarrow 2^{S_{i}}$\textit{\ is upper semicontinuous with
non-empty compact values;}

$G_{i}:\Omega _{i}\times S\rightarrow \mathcal{F}(S^{\prime })$ \textit{and} 
$g_{i}:S\rightarrow (0,1]$ \textit{are such that:}

\ \ (c) \textit{For each fixed }$(\omega _{i},y)\in \Omega _{i}\times S,$%
\textit{\ }$\lambda \rightarrow \{y\in S:\sup_{u\in (G_{i}(\omega
_{i},y))_{g_{i}(y)}}$Re$\langle u,\lambda _{i}(\omega _{i})-y\rangle >0\}:%
\mathcal{D}_{(X(\cdot ))_{z}}\rightarrow 2^{S}$\textit{\ is upper
semicontinuous with compact values;}

\ \ \ (\textit{d) for each fixed }$\omega _{i}\in \Omega _{i},$\textit{\ the
set}

\textit{\ }$\{\lambda \in \mathcal{D}_{(X(\cdot ))_{z}}:\sup_{y\in
(A_{i}(\omega _{i},\lambda ))_{a_{i}(\lambda )})}\sup_{u\in (G_{i}(\omega
_{i},y))_{g_{i}(y)}}$\textit{Re}$\langle u,\lambda _{i}(\omega
_{i})-y\rangle >0\}$\textit{\ is weakly open in }$\mathcal{D}_{(X(\cdot
))_{z}}$

\ \ \ (e) $\{\omega _{i}\in \Omega _{i}:\sup_{u\in (G_{i}(\omega
,y))_{g_{i}(y)}}\mathit{Re}\langle u,\lambda _{i}(\omega _{i})-y\rangle
>0\}\in \mathcal{Z}_{i}$ \textit{for each} $\lambda \in \mathcal{D}%
_{(X(\cdot ))_{z}}$\textit{.}

\textit{Then, there exists }$\lambda ^{\ast }\in \mathcal{D}_{(X(\cdot
))_{z}}$\textit{\ such that for every} $i\in I$:

i) $\lambda _{i}^{\ast }(\omega _{i})\in (A_{i}(\omega _{i},\lambda ^{\ast
}))_{a_{i}(\lambda ^{\ast })};$

ii) \textit{sup}$_{u\in (G_{i}(\omega _{i},y))_{g_{i}(y)}}Re\langle
u,\lambda _{i}^{\ast }(\omega _{i})-y\rangle \leq 0$\textit{\ for all }$y\in
(A_{i}(\omega _{i},\lambda ^{\ast }))_{a_{i}(\lambda ^{\ast })}$\textit{.}$%
\medskip $

\textit{Proof.} Let us define $\psi _{i}:\Omega _{i}\times \mathcal{D}%
_{(X(\cdot ))_{z}}\times S\rightarrow \mathbb{R}\cup \{-\infty ,+\infty \}$
by

$\psi _{i}(\omega _{i},\lambda ,y)=\sup_{u\in (G_{i}(\omega
_{i},y))_{g_{i}(y)}}$\textit{Re}$\langle u,\lambda _{i}(\omega
_{i})-y\rangle $ for each $(\omega _{i},\lambda ,y)\in \Omega _{i}\times 
\mathcal{D}_{(X(\cdot ))_{z}}\times S.$

We have that $\lambda _{i}(\omega _{i})\notin \{y\in S:\psi _{i}(\omega
_{i},\lambda ,y)>0\}$\textit{\ }for each fixed\textit{\ }$(\omega
_{i},\lambda )\in \Omega _{i}\times \mathcal{D}_{(X(\cdot ))_{z}}.$

All the hypotheses of Theorem 8 are satisfied. According to Theorem 8, there
exists $\lambda ^{\ast }\in \mathcal{D}_{(X(\cdot ))_{z}}$ such that $%
\lambda _{i}^{\ast }(\omega _{i})\in (A_{i}(\omega _{i},\lambda ^{\ast
}))_{a_{i}(\lambda ^{\ast })}$ for every $i\in I$ and

sup$_{y\in A_{i}(\omega _{i},\lambda ^{\ast }))_{a_{i}(\lambda ^{\ast
})}}\sup_{u\in (G_{i}(\omega _{i},y))_{g_{i}(y)}}\mathit{Re}\langle
u,\lambda _{i}^{\ast }(\omega _{i})-y\rangle \leq 0$ for every $i\in
I.\medskip $

If $|I|=1$, we obtain the following corollary.

\begin{corollary}
Let $S$ be a countable completely metrizable topological vector space and
let $(\Omega ,\mathcal{Z},\mathcal{\mu })$ be a measure space$.$\textit{\
Suppose that the following conditions are satisfied}$:$
\end{corollary}

\textit{a) the fuzzy mapping }$X:\Omega \rightarrow \mathcal{F}(S)$\textit{\
is measurable such that }$(X(\cdot ))_{z}:\Omega \rightarrow 2^{S}$ \textit{%
has compact values; }

\textit{b) the fuzzy mapping }$A$ \textit{is such that} \textit{for each} $%
\lambda \in \mathcal{D}_{(X(\cdot ))_{z}},$ $\mathit{(}A(\cdot ,\lambda
)_{a(\lambda )}:\Omega \rightarrow 2^{S}$\textit{\ is measurable\ and, for
all }$\omega \in \Omega ,$\textit{\ }$(A(\omega ,\cdot ))_{a(\cdot )}:%
\mathcal{D}_{(X(\cdot ))_{z}}\rightarrow 2^{S}$\textit{\ is upper
semicontinuous with non-empty compact values;}

$G:\Omega \times S\rightarrow \mathcal{F}(S^{\prime })$ \textit{and} $%
g:S\rightarrow (0,1]$ \textit{are such that:}

\ \ (c) \textit{For each fixed }$(\omega ,y)\in \Omega \times S,$\textit{\ }$%
\lambda \rightarrow \{y\in S:\sup_{u\in (G(\omega ,y))_{g(y)}}$Re$\langle
u,\lambda (\omega )-y\rangle >0\}:\mathcal{D}_{(X(\cdot ))_{z}}\rightarrow
2^{S}$\textit{\ is upper semicontinuous with compact values;}

\ \ \ (\textit{d) for each fixed }$\omega \in \Omega ,$\textit{\ the set}

\textit{\ }$\{\lambda \in \mathcal{D}_{(X(\cdot ))_{z}}:\sup_{y\in (A(\omega
,\lambda ))_{a(\lambda )})}\sup_{u\in (G(\omega ,y))_{g(y)}}$\textit{Re}$%
\langle u,\lambda (\omega )-y\rangle >0\}$\textit{\ is weakly open in }$%
\mathcal{D}_{(X(\cdot ))_{z}}$

\ \ \ (e) $\{(\omega ,\lambda ):\sup_{u\in (G(\omega ,y))_{g(y)}}\mathit{Re}%
\langle u,\lambda (\omega )-y\rangle >0\}\in \mathcal{Z}$ for each $\lambda
\in \mathcal{D}_{(X(\cdot ))_{z}}$\textit{.}

\textit{Then, there exists }$\lambda ^{\ast }\in \mathcal{D}_{(X(\cdot
))_{z}}$\textit{\ such that}:

i) $\lambda ^{\ast }(\omega )\in (A(\omega ,\lambda ^{\ast }))_{a(\lambda
^{\ast })};$

ii) \textit{sup}$_{u\in (G(\omega ,y))_{g(y)}}Re\langle u,\lambda ^{\ast
}(\omega )-y\rangle \leq 0$\textit{\ for all }$y\in (A(\omega ,\lambda
^{\ast }))_{a(\lambda ^{\ast })}$\textit{.}$\medskip $

We obtain the following random fixed point theorem as a particular case of
Theorem 8$.$

\begin{theorem}
Let $I$ be a non-empty finite set. For each $i\in I$, $S_{i}$ is a countable
complete metric space and $(\Omega _{i},\mathcal{Z}_{i})$ is a measurable
space. Let $(\Omega ,\mathcal{Z})$ be the product measurable space $\underset%
{i\in I}{(\prod }\Omega _{i},\underset{i\in I}{\tbigotimes }\mathcal{Z}_{i})$%
, and $\mu $ a probability measure on $(\Omega ,\mathcal{Z}).$\textit{\
Suppose that the following conditions are satisfied for each }$i\in I:$
\end{theorem}

\textit{a) the fuzzy mapping }$X_{i}:\Omega _{i}\rightarrow \mathcal{F}%
(S_{i})$\textit{\ is measurable and }$(X_{i}(\cdot ))_{z_{i}}:\Omega
_{i}\rightarrow 2^{S_{i}}$ \textit{has compact values; }

\textit{b) the fuzzy mapping }$A_{i}$ \textit{is such that} \textit{for each}
$\lambda \in \mathcal{D}_{(X(\cdot ))_{z}},$ $\mathit{(}A_{i}(\cdot ,\lambda
)_{a_{i}(\lambda )}:\Omega _{i}\rightarrow 2^{S_{i}}$\textit{\ is
measurable\ and, for all }$\omega _{i}\in \Omega _{i},$\textit{\ }$%
(A_{i}(\omega _{i},\cdot ))_{a_{i}(\cdot )}:\mathcal{D}_{(X(\cdot
))_{z}}\rightarrow 2^{S_{i}}$\textit{\ is upper semicontinuous with
non-empty compact values.}

\textit{Then, there exists }$\lambda ^{\ast }\in \mathcal{D}_{(X(\cdot
))_{z}}$\textit{\ such that for every} $i\in I$, $\lambda _{i}^{\ast
}(\omega _{i})\in (A_{i}(\omega _{i},\lambda ^{\ast }))_{a_{i}(\lambda
^{\ast })}$\textit{.}$\medskip $

If $|I|=1$, we obtain the following result.

\begin{theorem}
Let $S$ be a countable complete metric space and let $(\Omega ,\mathcal{Z},%
\mathcal{\mu })$ be a measure space. \textit{Suppose that the following
conditions are satisfied}$:$
\end{theorem}

\textit{a) the fuzzy mapping }$X:\Omega \rightarrow \mathcal{F}(S)$\textit{\
is measurable and }$(X(\cdot ))_{z}:\Omega \rightarrow 2^{S}$ \textit{has
compact values; }

\textit{b) the mapping }$A$ \textit{is such that} \textit{for each} $\lambda
\in \mathcal{D}_{(X(\cdot ))_{z}},$ $\mathit{(}A(\cdot ,\lambda )_{a(\lambda
)}:\Omega \rightarrow 2^{S}$\textit{\ is measurable\ and, for all }$\omega
\in \Omega ,$\textit{\ }$(A(\omega ,\cdot ))_{a(\cdot )}:\mathcal{D}%
_{(X(\cdot ))_{z}}\rightarrow 2^{S}$\textit{\ is upper semicontinuous with
non-empty compact values.}

\textit{Then, there exists }$\lambda ^{\ast }\in \mathcal{D}_{(X(\cdot
))_{z}}$\textit{\ such that} $\lambda ^{\ast }(\omega )\in (A(\omega
,\lambda ^{\ast }))_{a(\lambda ^{\ast })}$\textit{.}$\medskip $

\section{Classical systems of generalized quasi-variational inequalities and
random fixed point theorems}

If we consider classical corespondences in the last section, we obtain
several results concerning systems of generalized random quasi-variational
inequalities and random fixed points, which are new in the literature. We
present Theorem 14 as a consequence of Theorem 8.$\medskip $

\begin{theorem}
Let $I$ be a non-empty finite set. For each $i\in I$, $S_{i}$ is a countable
complete metric space and $(\Omega _{i},\mathcal{Z}_{i})$ is a measurable
space. Let $(\Omega ,\mathcal{Z})$ be the product measurable space $\underset%
{i\in I}{(\prod }\Omega _{i},\underset{i\in I}{\tbigotimes }\mathcal{Z}_{i})$%
, and let $\mu $ be a probability measure on $(\Omega ,\mathcal{Z}).$\textit{%
\ Suppose that the following conditions are satisfied for each }$i\in I:$
\end{theorem}

\textit{a) the correspondence }$X_{i}:\Omega _{i}\rightarrow 2^{S_{i}}$%
\textit{\ is measurable and }$(X_{i}(\cdot ))_{z_{i}}:\Omega _{i}\rightarrow
2^{S_{i}}$ \textit{has compact values; }

\textit{b) the correspondence }$A_{i}$ \textit{is such that} \textit{for each%
} $\lambda \in \mathcal{D}_{X},$ $A_{i}(\cdot ,\lambda ):\Omega
_{i}\rightarrow 2^{S_{i}}$\textit{\ is measurable\ and, for all }$\omega
_{i}\in \Omega _{i},$\textit{\ }$A_{i}(\omega _{i},\cdot ):\mathcal{D}%
_{X}\rightarrow 2^{S_{i}}$\textit{\ is upper semicontinuous with non-empty
compact values;}

Let us \textit{assume that the mapping }$\psi _{i}:\Omega _{i}\times 
\mathcal{D}_{X}\times S\rightarrow \mathbb{R}\cup \{-\infty ,+\infty \}$%
\textit{\ is such that:}

\ \ \ (\textit{c) }$\lambda \rightarrow \{y\in Y:\psi _{i}(\omega ,\lambda
,y)>0\}:\mathcal{D}_{X}\rightarrow 2^{S}$\textit{\ is upper semicontinuous
with compact values on }$\mathcal{D}_{X}$\textit{\ for each fixed }$\omega
_{i}\in \Omega _{i};$

\ \ \ (\textit{d) }$\lambda _{i}(\omega _{i})\notin \{y\in S:\psi
_{i}(\omega _{i},\lambda ,y)>0\}$\textit{\ for each fixed }$(\omega
_{i},\lambda )\in \Omega _{i}\times \mathcal{D}_{X};$

\ \ \ (\textit{e) for each }$\omega _{i}\in \Omega _{i},$\textit{\ }$%
\{\lambda \in \mathcal{D}_{X}:\alpha _{i}(\omega _{i},\lambda )>0\}$\textit{%
\ is weakly open in} $\mathcal{D}_{X},$\textit{\ where }$\alpha _{i}:\Omega
_{i}\times \mathcal{D}_{X}\rightarrow \mathbb{R}$\textit{\ is defined by }$%
\alpha _{i}(\omega _{i},\lambda )=\sup_{y\in A_{i}(\omega _{i},\lambda
)}\psi _{i}(\omega _{i},\lambda ,y)$\textit{\ for each }$(\omega
_{i},\lambda )\in \Omega _{i}\times \mathcal{D}_{X};$

\ \ \ (\textit{f) }$\{\omega _{i}:\alpha _{i}(\omega _{i},\lambda )>0\}\in 
\mathcal{Z}_{i}$ for each $\lambda \in \mathcal{D}_{X}$\textit{.}

\textit{Then, there exists }$\lambda ^{\ast }\in \mathcal{D}_{X}$\textit{\
such that for every} $i\in I$,

i) $\lambda _{i}^{\ast }(\omega _{i})\in A_{i}(\omega _{i},\lambda ^{\ast
}); $

\textit{ii) sup}$_{y\in A_{i}(\omega ,\lambda ^{\ast })}\psi _{i}(\omega
_{i},\lambda ^{\ast },y)\leq 0$\textit{.}$\medskip $

The next theorems concern the existence of the random fixed point for
correspondences with values in complete countable metric spaces.$\medskip $

\begin{theorem}
Let $I$ be a non-empty finite set. For each $i\in I$, $S_{i}$ is a countable
complete metric space and $(\Omega _{i},\mathcal{Z}_{i})$ is a measurable
space. Let $(\Omega ,\mathcal{Z})$ be the product measurable space $\underset%
{i\in I}{(\prod }\Omega _{i},\underset{i\in I}{\tbigotimes }\mathcal{Z}_{i})$%
, and let $\mu $ be a probability measure on $(\Omega ,\mathcal{Z}).$\textit{%
\ Suppose that the following conditions are satisfied for each }$i\in I:$
\end{theorem}

\textit{a) the correspondence }$X_{i}:\Omega _{i}\rightarrow 2^{S_{i}}$%
\textit{\ is measurable with compact values; }

\textit{b) for each} $\lambda \in \mathcal{D}_{X},$ $A_{i}(\cdot ,\lambda
):\Omega _{i}\rightarrow 2^{S_{i}}$\textit{\ is measurable\ and, for all }$%
\omega _{i}\in \Omega _{i},$\textit{\ }$A_{i}(\omega _{i},\cdot ):\mathcal{D}%
_{X}\rightarrow 2^{S_{i}}$\textit{\ is upper semicontinuous with non-empty
compact values.}

\textit{Then, there exists }$\lambda ^{\ast }\in \mathcal{D}_{X}$\textit{\
such that for every} $i\in I$, $\lambda _{i}^{\ast }(\omega _{i})\in
A_{i}(\omega _{i},\lambda ^{\ast })$\textit{.}$\medskip $

If $|I|=1$, we obtain the following result.

\begin{theorem}
Let $S$ be a countable complete metric space and let $(\Omega ,\mathcal{Z},%
\mathcal{\mu })$ be a measure space. \textit{Suppose that the following
conditions are satisfied}$:$
\end{theorem}

\textit{a) the correspondence }$X:\Omega \rightarrow 2^{S}$\textit{\ is
measurable with compact values; }

\textit{b) for each} $\lambda \in \mathcal{D}_{X},$ $A(\cdot ,\lambda
):\Omega \rightarrow 2^{S}$\textit{\ is measurable\ and, for all }$\omega
\in \Omega ,$\textit{\ }$A(\omega ,\cdot ):\mathcal{D}_{X}\rightarrow 2^{S}$%
\textit{\ is upper semicontinuous with non-empty compact values.}

\textit{Then, there exists }$\lambda ^{\ast }\in \mathcal{D}_{X}$\textit{\
such that} $\lambda ^{\ast }(\omega )\in A(\omega ,\lambda ^{\ast })$\textit{%
.}$\medskip $

\begin{center}
\bigskip
\end{center}

\end{document}